\documentclass[a4paper,10pt]{article}
\usepackage{amssymb}
\usepackage{amsmath,amsfonts,amsthm,amssymb}
\usepackage[dvips]{graphics}
\usepackage{epsfig}
\usepackage{indentfirst}
\usepackage{cancel,soul}

\usepackage{color}
\usepackage{bm}
\allowdisplaybreaks

\usepackage[colorlinks=true,citecolor=red,linkcolor=blue,%
urlcolor=RubineRed,pdfpagetransition=Blinds,pdftoolbar=false,pdfmenubar=false]{hyperref}

\newtheoremstyle{theorem}
  {10pt}          
  {10pt}  
  {\sl}  
 {}
  {\bf}  
  {. }    
  { }    
  {}     
\theoremstyle{theorem}

\newtheorem{theorem}{Theorem}[section]

 \newtheorem{lemma}{Lemma}[section]
 
 \newtheorem{remark}{Remark}[section]
  
\numberwithin{equation}{section}

\newtheoremstyle{defi}
  {10pt}          
  {10pt}  
  {\rm}  
  {}  
  {\bf}  
  {. }    
  { }    
  {}     
\theoremstyle{defi}

\def\x{{\textbf {\textit x}}}

\def\vi{{\textbf {\textit v}}}

\begin{document}
\baselineskip = 16pt

\title{\bf Energy equality for the compressible Primitive Equations with vacuum}

\author{\v{S}\'{a}rka Ne\v{c}asov\'{a}$^1$
\footnote{Email: matus@math.cas.cz} \ \ \
Mar\'{\i}a \'{A}ngeles Rodr\'{\i}guez-Bellido$^{2}$\footnote{Email: angeles@us.es}\\
Tong Tang$^{3}$ \footnote{Email: tt0507010156@126.com}\\
{\small  1. Institute of Mathematics of the Academy of Sciences of the Czech Republic,} \\
{\small \v Zitn\' a 25, 11567, Praha 1, Czech Republic}\\
{\small 2. Dpto. Ecuaciones Diferenciales y An\'alisis Num\'erico and IMUS, Universidad de Sevilla}\\
{\small Campus de Reina Mercedes, Facultad de Matem\'aticas, C/ Tarfia, s/n- 41012 Sevilla, Spain}\\
{\small 3. School of Mathematical Science}\\
{\small Yangzhou University, Yangzhou 225002, P.R. China}\\
\date{}}

\maketitle

\begin{abstract}
The paper deals with the problem of the energy conservation for the weak solutions to the compressible Primitive Equations (CPE) system with degenerate viscosity.
The sufficient conditions on the regularity of weak solutions for the energy equality are obtained
even for the case when the solutions may include vacuum. In this paper, we show two theorems, the first one gives regularity in the classical isotropic Sobolev and Besov spaces. The second one  states regularity in the anisotropic spaces. We obtain new regularity results in the second theorem due to the special structure of CPE system, which are in  contrast to compressible Navier-Stokes equations.

\vspace{0.5cm}


\medskip

\end{abstract}

\maketitle
\section{Introduction}\setcounter{equation}{0}
In the paper we focus on  the relationship between regularity and conservation of energy
for the compressible Primitive Equations (CPE) system in the periodic domain $\mathbb{T}^3$. The system is described the following form
\begin{equation}\label{1.1}
\left\{
\begin{array}{llll}
\partial_{t}\rho+\text{div}_{\bf x}(\rho \mathbf{v})+\partial_z(\rho w)=0, \\
\partial_t(\rho \mathbf{v})+\textrm{div}_{\bf x}(\rho\mathbf{v}\otimes\mathbf{v})+\partial_z(\rho\mathbf vw)+\nabla_{\bf x} p(\rho)={\rm div}_{\bf x}(\mu(\rho)\nabla_{\bf x}\mathbf v)+\partial_{z}(\mu(\rho)\partial_z\mathbf v),\\
\partial_zp(\rho)=0,
\end{array}\right.
\end{equation}
where $\rho, \mathbf u, p$ represent the density, velocity and pressure, respectively.
The velocity is defined as $\mathbf{u}=(\mathbf v,w)$, $\mathbf v(t,{\bf x},z)\in\mathbb{R}^2$
and $w(t,{\bf x},z)\in\mathbb{R}$ represent the horizonal velocity and vertical velocity respectively and
where ${{\bf x} \in \mathbb{R}^2}$ denotes the horizontal direction and $z$ denotes the vertical direction. Here we assume the linear dependency of the viscosity on the density:  $\mu(\rho)=\rho$. From $(\ref{1.1})_3$, we can assume  that the density is independent of $z$, which means $\rho=\rho(t,{\bf x})$.
Such assumption can be found also in \cite{er1,ga}, where the authors use a change of variables to derive the density $\rho$ is independent of $z$.
Moreover we consider pressure as $p(\rho)=\rho^\gamma\hspace{5pt}(\gamma>1)$, and the density could vanish ($\rho\ge 0$).

{\bf Primitive Equations (PE)}\footnote{PE we mean incompressible primitive equations.}
system is an important model which is widely used in the geophysical research to describe and analyze the phenomena of atmosphere and ocean. It is derived from the Navier-Stokes or Euler systems by asymptotic analysis or numerical approximation. Let us briefly recall some important results for such a system.

During the last decades, there is a large amounts of literature concerning the  rigourous mathematical justification of deriving PE model.
More precisely, for incompressible (PE) system, Az\'{e}rad and Guill\'{e}n \cite{a} proved that the incompressible Navier-Stokes equations converge to PE in the sense of weak solutions.
 Further, Li and Titi \cite{li} showed the convergence of the weak solutions of incompressible Navier-Stokes equations  to the strong solutions of PE. Based on \cite{a,li},
 Donatelli and Juhasz \cite{do} gave a justification that PE model with the pollution effect is the hydrostatic limit of the Navier-Stokes equations with an advection-diffusion equation.
 Grenier \cite{gr} used the energy estimates and Brenier \cite{by} used the relative entropy inequality to prove that the smooth solutions of incompressible Euler system converge to smooth solutions of inviscid PE. Precisely, Brenier \cite{b} proved the existence of smooth solutions in two-dimensions under the convex horizontal velocity assumptions.  Later, Masmoudi and Wong \cite{m} extended Brenier's result, removing the convex horizontal velocity assumptions.

 On the other hand, Ersoy et al. \cite{er1} used the asymptotic analysis to deduce the {\it compressible primitive system (CPE)} with degenerate viscosity coefficients. Gao, Ne\v casov\' a and Tang \cite{gao} deduced the CPE from anisotropic Navier-Stokes equations with constant viscosity coefficient.
\smallskip

One of the typical features of the PE model is that there is no information for the vertical velocity in the momentum equation and the vertical velocity is determined by the horizontal velocity through the continuity equation that here reduces to incompressible constraint.
Therefore, the mathematical and numerical study of the PE model were unsolved until 1990s when Lions, Teman and Wang \cite{l1,l2} established fundamental results in this field.
After that{\color{red},} Guill\'en-Gonz\'alez, Masmoudi and Rodr\'{\i}guez-Bellido \cite{gu} proved the local existence of strong solutions and uniqueness with some interesting anisotropic estimates.
The celebrated result was made by Cao and Titi \cite{c1}, where they first proved the global well-posedness of PE in the three dimensional case. Then, by virtue of semigroup method, Hieber and Kashiwabara \cite{h} extended this result relaxing the smoothness on the initial data.
\smallskip

{\bf Compressible primitive equations} system (CPE) has been studied more recently. Gatapov and Kazhikhov \cite{ga} proved the global existence of weak solutions with the constant viscosity coefficients in 2D case.
Liu and Titi \cite{liu1,liu3} proved the local existence of strong solutions in 3D case and considered the zero Mach number limit of CPE.
Regarding the degenerate viscosity case, Ersoy et al. \cite{er1}, and
Tang and Gao \cite{tang} showed the stability of weak solution.
The stability means that a subsequence of weak solutions will converge to another weak solutions if it satisfies some uniform bounds.
Liu and Titi \cite{liu2} and Wang et al. \cite{w}, independently, used the B-D entropy to prove the global existence of weak solutions.
Readers can refer to Bresch et al.~\cite{bg}, Cao et al.~\cite{c4}, Li and Titi~\cite{li}, the book of Temam and Ziane {\cite{t} and references therein for more physical background and other interesting mathematical results.
\smallskip

The existence of weak solutions is a fundamental question in PDE, especially for fluids models. The case of the incompressible fluids was investigated already in thirties by famous works of Leray \cite{leray}. He proved the global existence of weak solutions to incompressible Navier-Stokes equations. On the other hand, the proof of existence for the compressible case is going back to the nineties by Lions \cite{pl2} and Feireisl et al. \cite{e1,e2,e3} who proved the existence of global weak solutions of compressible isentropic case and later to full system.

{\it Generally, it is not  known whether weak solutions satisfy the principle of conservation of energy for both incompressible and compressible fluids.}
It is a nature question how high regularity for weak solutions is needed to obtain the energy equality.
\smallskip

 Onsager \cite{o} gave a famous conjecture that the three dimensional {\it incompressible Euler equations} conserve energy if the velocity
$\mathbf u\in L^3((0,T);C^{0,\alpha}(\mathbb{T}^3))$ with $\alpha>\frac{1}{3}$.
 The second part of conjecture said that there exist weak solutions of the Euler equation for $\alpha<\frac{1}{3}$ which do not conserve energy.

The second part of Onsager's conjecture has been underlined by the celebrated work of Scheffer \cite{sc} and Shnirelman \cite{sh}.
And a series of breakthrough papers were done by De Lellis and Sz\'{e}kelyhidi \cite{de1,de2,de3} by virtue of convex integration. Recently, this part is fully solved by Isett \cite{i}, Buckmaster et. al. \cite{bu}.

The first part of  Onsager's conjecture was expressed as energy conservation, which is a lively direction of research at the contemporary mathematical society. It was proved by Constantin, E and Titi in \cite{co} (also by  Eyink \cite{ey}, and the work of Duchon and Robert \cite{du}), stating that if
${\bf u}$ belongs to $L^3([0,T];B^{\alpha,\infty}_3(\mathbb{T}^3))\cap C([0,T];L^2(\mathbb{T}^3))$ with $\alpha>\frac{1}{3}$,
then the energy is conserved. Cheskidov et al. \cite{ch}, and Fjordholm and Wiedemann \cite{f} made the sharpest result in optimal Besov spaces.
Their main idea is using the suitable commutator estimates for incompressible Euler system.
These results were extended to the bounded domain by Bardos, Titi and Wiedemann \cite{b1}, Drivas and Nguyen \cite{dr}.
Later this result was extended to incompressible {\it inhomogeneous} Euler equations by Feireisl et al. \cite{e4}.
Precisely speaking, the authors in \cite{e4} mollified the weak solutions of the density $\rho$ and velocity $\mathbf u$
and stated that if ${\bf u}$ belongs to $B^{\alpha,\infty}_p([0,T];\times \mathbb{T}^3)$,$\varrho, \varrho {\bf u} \in B^{\beta,\infty}_q([0,T];\times \mathbb{T}^3)$, $p \in L^{p^{*}}_{loc}((0,T) \times \mathbb{T}^3)$ for some $1\leq p, q \leq \infty$, $p^*$ is the conjugate of $p$, and $ 0\leq \alpha, \beta \leq 1$,
then the energy is locally conserved.
On the other hand, Leslie and Shvydkoy \cite{le} extended these results into the inhomogeneous Navier-Stokes case.
Moreover, Gwiazda et al. \cite{gw} proved the corresponding result to general system of first order conservation laws.
\smallskip

Concerning {\it compressible Euler system} such result goes to Feireisl et al \cite{e4}. Inspired by Constantin et al. \cite{co}, Feireisl et al. \cite{e4} proved the energy conservation for compressible Euler system with initial data containing vacuum.
Yu \cite{y} used the Lions's commutator estimates to show energy conservation for compressible Navier-Stokes equations with degenerate viscosity but without vacuum.
Nguyen et al. \cite{n} extended Yu's result with weaker regularity condition in bounded domain.
Recently, Akramov et al. \cite{ak} proved the corresponding results for compressible Euler and Navier-Stokes system by removing the key assumptions that the pressure is a twice continuously differentiable function.
One can refer to \cite{b2,bo1,de} for other related works.
\smallskip

Comparing with fruitful results for either incompressible and compressible Euler
 and Navier-Stokes system,
 there are a few results about PE model due to its special structure.
There are many differences at the mathematical structure between Navier-Stokes equations and PE model, due to the hydrostatic approximation.
There is {\bf no information for the vertical velocity in the momentum equation of PE model}, so it is very difficult to analyze the PE model.
As far as we know, there is only one result due to Boutros, Markfelder and Titi \cite{bo}, in which they showed the energy conservation for {\bf the incompressible inviscid PE}.
They gave three types of weak solutions, where their vertical velocity is determined by horizontal velocity by virtue of incompressible conditions.
One of the sufficient condition for energy conservation is horizonal velocity $\mathbf v\in L^4((0,T);B^{\alpha,\infty}_{4}(\mathbb{T}^3))$ ($\alpha>\frac{1}{2}$) and vertical velocity $w\in L^2((0,T);L^2(\mathbb{T}^3))$.
It is interesting to find that the index for velocity regularity is different from conventional Onsager's exponent, which coincides with what stated in \cite{bo1}, ``Several of these exponents were different from $\frac{1}{3}$, which is the Onsager exponent for the Euler equations.
This means that the Onsager exponent is not universal and is plausibly determined by the
regularity of the coefficient functions of the nonlinearity..." It also shows the PE has different mathematical structure and corresponding different results compared with Navier-Stokes equations.
\smallskip


{\bf There is no result concerning the energy equality for CPE system.}
The goal of this work is to investigate the energy equality for CPE system  under some regularity conditions.
{\bf There are some distinct differences between our work and Feireisl et al. \cite{e4},
and  Akramov  et al. \cite{ak}.}
Theorem \ref{t2.1} seems similar and parallel to the results in \cite{e4,ak}. However, in contrast with the compressible Euler and Navier-Stokes system, we just have regularity of the horizontal velocity in CPE. The only way to obtain the information for the vertical velocity is through the continuity equation instead of the incompressible constraint. It means it is a hurdle to control the estimate of vertical velocity in the isotropic space. In order to overcome this difficulty, we need the additional assumption $\rho\mathbf v\in L^3((0,T);H^1)$ in Theorem \ref{t2.1}. Fortunately, by the deep investigation of the structure of CPE, we are able to remove this restriction by virtue of the anisotropic Besov space regularity in Theorem \ref{t2.2}. Therefore, we can obtain the result on the anisotropic Besov space and solve the corresponding problems.

The paper is organized as follows:
In Section \ref{S2}, we introduce the Besov space, give some useful lemmas and state the main theorems.
Section \ref{S3} is devoted to the proof of the energy conservation.
\vskip 0.5cm

\section{Preliminaries} \label{S2}

\subsection{Function spaces}
Let $\Omega:=(0,T)\times\mathbb{T}^3$, and we define the Besov space:
$B^{\alpha,\infty}_{p}(\Omega)$ with $1\leq p<\infty$, $0<\alpha<1$,
normed by
\begin{align*}
\|w\|_{B^{\alpha,\infty}_{p}(\Omega)}
=
\|w\|_{L^p(\Omega)}+\sup_{{\boldsymbol \xi}\in\Omega}\{|{\boldsymbol \xi}|^{-\alpha}\|w(\cdot+{\boldsymbol \xi})-w\|_{L^p(\Omega\cap(\Omega-{\boldsymbol \xi}))}\}
\end{align*}
is finite. In order to understand the anisotropic structure in CPE, inspired by \cite{bo,gu}, we introduce the following anisotropic spaces:
\begin{align}\label{2.1}
\|w\|_{H^1_{\bf x}(\mathbb{T}^2)L^2_z(\mathbb{T}^1)}=\|\|w( {\bf x},\cdot)\|_{L^2(\mathbb{T}^1)}\|_{H^1_{\bf x}(\mathbb{T}^2)}
\end{align}
and
considering ${\boldsymbol \xi}=(\tau,{\boldsymbol \xi}_h,\xi_z)$,

\begin{align}\label{2.2}
\|w\|_{B^{\alpha,\infty}_{p}((0,T);B^{\alpha,\infty}_{p}(\mathbb{T};B^{\beta,\infty}_{p}(\mathbb{T}^2))}
=&\|w\|_{L^p(\Omega)}+\sup_{{\tau}\in \mathbb{R}\setminus\{0\}}\{|{\tau}|^{-\alpha}\|w(\cdot+(\tau,0,0))-w\|_{L^p(\Omega)\cap(\Omega-(\tau,0,0))}\}\nonumber\\
&+ \sup_{{\boldsymbol\xi}\in \mathbb{R}^3\setminus\{0\}}\{|{\boldsymbol \xi}|^{-\alpha}\|w(\cdot+(0,{\boldsymbol \xi}_h,\xi_z))-w\|_{L^p(\Omega)\cap(\Omega-(0,{\boldsymbol \xi}_h,\xi_z))}\}\nonumber\\
&+\sup_{{\boldsymbol \xi}_h\in\mathbb{R}^2\setminus\{0\}}\{
|{\boldsymbol \xi}|^{-\beta}
\|w(\cdot+(0,{\boldsymbol \xi}_h,0))-w\|_{L^p(\Omega)\cap(\Omega-(0,{\boldsymbol \xi}_h,0 ))}\},
\end{align}
where $0<\alpha<\beta<1$.

The space $H^1_{\bf x}(\mathbb{T}^2)L^2_z(\mathbb{T}^1)$ (which appears in Theorem \ref{t2.2}) means its horizontal regularity is $H^1$ and vertical regularity is $L^2$.

Let $\eta\in C^\infty_c(\mathbb{R}^{d+1})$ ($d$ is the number of the space dimension) be a standard mollification kernel and set
\begin{align*}
\eta^{\varepsilon}({\bf y})=\frac{1}{\varepsilon^{d+1}}\eta\left(\frac{ {\bf y}}{\varepsilon}\right),\hspace{3pt} w^{\epsilon}=\eta^\varepsilon\ast w,\hspace{3pt}
f^\varepsilon(w)=f(w)\ast\eta^\varepsilon,
\end{align*}
where $w^{\varepsilon}$ is well-defined on $\Omega^{\varepsilon}=\{{\bf y} \in \Omega: \, d({\bf y},\partial\Omega)>\varepsilon\}$. Now we recall some useful lemmas which will be frequently used throughout the paper.

\bigskip

\begin{lemma}(\cite{e4,de})\label{2.4}
For any function $u\in B^{\alpha,\infty}_{p}(\Omega)$, we have
\begin{align} 
&\|u(\cdot+{\boldsymbol \xi})-u\|_{L^p(\Omega\cap(\Omega-{\boldsymbol \xi}))}\leq|{\boldsymbol \xi}|^\alpha\|u\|_{B^{\alpha,\infty}_{p}(\Omega)}, \nonumber\\
&\|u^\varepsilon-u\|_{L^p(\Omega)}\leq\varepsilon^\alpha\|u\|_{B^{\alpha,\infty}_{p}(\Omega)}, \nonumber\\
&\|\nabla u^\varepsilon\|_{L^p(\Omega)}\leq\varepsilon^{\alpha-1}\|u\|_{B^{\alpha,\infty}_{p}(\Omega)}, \nonumber
\end{align}
where $\nabla$ stands for space-time gradient.
\end{lemma}

\begin{remark}\label{r2.2}
It is easy to prove the following equality \cite{e4}
\begin{align*}
&f^\varepsilon g^\varepsilon
-(fg)^\varepsilon=
(f^\varepsilon-f)(g^\varepsilon-g)\\
&
-\displaystyle\int^\varepsilon_{-\varepsilon}
\int_{\mathbb{T}^3}
\- \- \eta^\varepsilon(\tau,{\boldsymbol \xi}_h,\xi_z)
\left(
f(t-\tau,{\bf x}-{\boldsymbol \xi}_h,z-\xi_z)-
f (t,{\bf x},z)
\right)
\left(
g(t-\tau,{\bf x}-{\boldsymbol \xi}_h,z-\xi_z)-g(t,{\bf x},z)
\right)
d{\boldsymbol \xi}_h  d\xi_z   d\tau.
\end{align*}

\end{remark}

\begin{lemma} (\cite{ba} Besov space embedding) \label{l2.2}
Let $s\in\mathbb{R}$, $1\leq p,r\leq\infty$, then
\begin{align} 
&B^{s,r}_p\hookrightarrow B^{s-d(\frac{1}{p}-\frac{1}{p_1}),r}_{p_1},\hspace{3pt}p_1\geq p.
\nonumber
\end{align}
\end{lemma}
It is easy to obtain the following remark:
\begin{remark}\label{r2.3}
$H^1(\mathbb{T}^3)\hookrightarrow B^{1,\infty}_{2}(\mathbb{T}^3)\hookrightarrow B^{\alpha,\infty}_{3}(\mathbb{T}^3)$\hspace{3pt}($\alpha>\frac{1}{2}$).
\end{remark}

The following two lemmas can be proved as \cite{de}.
\begin{lemma}\label{l2.3}
Let $u\in B^{\alpha,\infty}_{p}(\Omega,\mathbb{R}^m)$. Suppose $f:\mathbb{R}^m\rightarrow\mathbb{R}$ is a $C^1$ function with $\frac{\partial f}{\partial u_i}\in L^\infty$ $i=1...m$. Then
\begin{align}
\|\nabla f(u^\varepsilon)\|_{L^p(\Omega)}\leq C\varepsilon^{\alpha-1}\|u\|_{B^{\alpha,\infty}_{p}(\Omega)},
\nonumber
\end{align}
where $\nabla$ stands for space-time gradient
\end{lemma}

\begin{lemma}\label{l2.4a}
Let $1\leq q<\infty$ and suppose $v\in L^{2q}((0,T)\times\mathbb{T}^d;\mathbb{R}^k)$ and $f\in C^2(\mathbb{R}^k;\mathbb{R}^N)$. If
\begin{align*}
\sup_{i,j}\|\frac{\partial^2f}{\partial v_i\partial v_j}\|_{L^\infty}<\infty,
\end{align*}
then there exists a constant $C>0$ such that
\begin{align}
\|f(v^\varepsilon)-f^\varepsilon(v)\|_{L^q}<C\big{(}\|v^\varepsilon-v\|^2_{L^{2q}}
+
\sup_{(s,{\bf y})\in {\rm supp}(\eta_\varepsilon)}\|v(\cdot,\cdot)-v(\cdot-s,\cdot-{\bf y})\|^2_{L^{2q}}\big{)}.
\nonumber
\end{align}
\end{lemma}

The following lemma, known as  Lions's commutator lemma (\cite{y} Lemma 2.3 and \cite{la} Lemma 3.2), will be very useful for the estimates below.

\begin{lemma}\label{l2.5}
Let $\partial$ be a partial derivative in space or time. Let $f$, $\partial f\in L^p(\mathbb{R}^{+}\times\Omega)$,
$g\in L^q(\mathbb{R}^{+}\times\Omega)$ with $1 \leq p,q\leq\infty$, and $\frac{1}{p}+\frac{1}{q}\leq1$. Then, we have
\begin{align}
\|\partial(fg)^\varepsilon-\partial(fg^\varepsilon)\|_{L^r}\leq C\|\partial f\|_{L^{p}} \|g\|_{L^{q}}
\nonumber
\end{align}
for some constant $C>0$ independent of $\varepsilon$, $f$ and $g$, with $\frac{1}{r}=\frac{1}{p}+\frac{1}{q}$. In addition,
\begin{align}
\partial(fg)^\varepsilon-\partial(fg^\varepsilon)\rightarrow0,\hspace{3pt}\text{in}\hspace{2pt} L^r,
\nonumber
\end{align}
as $\varepsilon\rightarrow0$, if $r<\infty$.
\end{lemma}
\begin{remark}
The original version in \cite{y} is $\|(\partial(fg))^\varepsilon-\partial(fg^\varepsilon)\|_{L^r}\leq C\|\partial f\|_{L^{p}} \|g\|_{L^{q}}$.
Using Lemma 3.1 in \cite{la}, $\partial(fg)^\varepsilon=(\partial(fg))^\varepsilon$, we get the statement of Lemma \ref{l2.5} as $$\|\partial(fg)^\varepsilon-\partial(fg^\varepsilon)\|_{L^r}=\|(\partial(fg))^\varepsilon-\partial(fg^\varepsilon)\|_{L^r}.$$
\end{remark}

\begin{remark}
The convergence of commutators is a crucial step in the proof of the energy conservation. The different assumptions of the regularity of weak solutions require various methods of proof.
For example, there are big differences in the proof between \cite{y} and \cite{ak,e4,gh}. Yu \cite{y} adopts Lions's commutator lemma (Lemma \ref{l2.5}) to show the convergence
to zero of all commutators with the regularity of weak solution in Sobolev space.
 Akramov et al. \cite{ak,gh} and Feireisl et al. \cite{e4} consider regularity in Besov space and utilize the decay property of Besov space (Lemma \ref{2.4}) to prove the corresponding convergence. Due to the special structure of CPE, where the vertical velocity $w$ does not satisfied a momentum equation, although these previous ideas are used, the control of the horizontal and vertical viscosity commutators $R^\varepsilon_4$ and $R^\varepsilon_5$ (and also the terms involving $w$) is made differently.
\end{remark}


\subsection{Main result}
Our main results can be stated as follows:

\begin{theorem}\label{t2.1}
Let $(\rho,\mathbf u)$ be a solution of CPE \eqref{1.1} in the sense of distribution.
Let us assume $\mathbf v\in B^{\alpha,\infty}_3((0,T);B^{\alpha,\infty}_{3}(\mathbb{T}^3))$,
$\rho,\rho\mathbf v\in B^{\beta,\infty}_3((0,T);B^{\beta,\infty}_{3}(\mathbb{T}^3)$,
and
\begin{align*}
0\leq\underline{\rho}\leq\rho\leq\overline{\rho},\hspace{3pt}\mbox{\text a.e.   in}\hspace{2pt}(0,T)\times\mathbb{T}^3
\end{align*}
for some constants $\underline{\rho},\overline{\rho}$ and $0\leq\alpha,\beta\leq1$, and $\alpha>\frac{5}{6}$.
Moreover, we assume the additional regularity $w\in L^2(0,T;L^2(\mathbb{T}^3))$, $\mathbf v\in L^{2}((0,T);H^1(\mathbb{T}^3))$, and $p\in C[\underline{\rho},\overline{\rho}]$.
Then the energy is locally conserved i.e.
\begin{equation}\label{energy1}
\begin{array}{r}
\displaystyle\partial_t
\left(
\frac{1}{2}\rho|\mathbf v|^2+P(\rho)
\right)+\rho|\nabla_{\bf x}\mathbf v|^2+\rho|\partial_z\mathbf v|^2
+{\rm div}_{\bf x} \left[
(\frac{1}{2}\rho|\mathbf v|^2+p(\rho)+P(\rho))\mathbf v
-\rho(\nabla_{\bf x}\mathbf v)\mathbf v
\right]
\\
\noalign{\vspace{-1ex}}\\
\displaystyle
+
\partial_z
\left[
\left(
\frac{1}{2}\rho|\mathbf v|^2+p(\rho)+P(\rho)
\right)
w
-
\rho\partial_z\mathbf v\cdot\mathbf v
\right]
=0,
\end{array}
\end{equation}
in the sense of distributions on $(0,T)\times\mathbb{T}^3$. Note that a new term $P(\rho)$ is introduced, defined as $P(\rho)=\rho\displaystyle\int^\rho_1\frac{p(s)}{s^2}ds$ and called the pressure potential.
\end{theorem}

\begin{remark}\label{t2.5}
If $\mathbf v\in B^{\alpha,\infty}_4((0,T);B^{\alpha,\infty}_{4}(\mathbb{T}^3))$ is assumed in the proof of Theorem \ref{t2.1}, we could still get the same result,
but with the relaxed index $\alpha>\frac{1}{2}$. It is similar to Titi et al. \cite{bo} , where they assume $\mathbf u\in L^4((0,T);B^{\alpha,\infty}_{4}(\mathbb{T}^3)))$  $(\alpha>\frac{1}{2})$ at Proposition 3.1.
We will show the key point in Remark \ref{t3.2} in the proof of Theorem \ref{t2.1}.
\end{remark}

\begin{remark}
In \cite{bo} Proposition 5.1, Titi et al. assume 
$\mathbf v\in L^3((0,T);W^{1,p}(\mathbb{T}^3))$ $(p>6)$, as they deal with the inviscid PE system. Inspired by methods and techniques introduced from Akramov et al. \cite{ak}, we can replace previous assumption on the regularity of horizontal velocity by $\mathbf v\in L^2((0,T);H^1(\mathbb{T}^3))$ in Theorem \ref{t2.1}.
\end{remark}

\begin{remark}\label{l2.4}
By virtue of the Bresch-Desjardins entropy, Liu and Titi \cite{liu2} obtained the extra estimates and the existence of global weak solutions to CPE with degenerate viscosity.
However, the existence of weak solutions to CPE with constant viscosity is still an open problem.
\end{remark}

\begin{remark}{ \bf [Differences with previous results]}
Vacuum is an important issue and challenge in the study of compressible Euler and Navier-Stokes system. Titi et al. \cite{bo} consider the energy conservation for the incompressible PE. Therefore, there is no vacuum difficulty in \cite{bo}. However, we consider the CPE with density, concerning the issue of vacuum with the assumption \begin{align*}
0\leq\underline{\rho}\leq\rho\leq\overline{\rho},\hspace{3pt}{\text a.e. in}\hspace{2pt}(0,T)\times\mathbb{T}^3.
\end{align*}
We should claim that we can relax the temporal regularity, if we exclude the vacuum, that is $\underline{\rho}>0$. In this case, the test function used in the momentum equation should be $\frac{(\rho \mathbf v)^\epsilon}{\rho^\epsilon}$ instead of $\mathbf v^\epsilon$, which is stated clearly in \cite{e4} Remark 4.2 and the introduction in \cite{n}.

Moreover, in contrast to the incompressible PE system, the pressure $p=p(\rho)$ is no longer a Lagrange multiplier, but a function of the density and needs to be controlled.
Last but not least, we have that the dissipative term and the viscosity depends on the density, which brings the new commutator estimates.
Therefore, the conventional method used in incompressible PE \cite{bo} can not be applied straightforwardly to the CPE model.
\end{remark}

Inspired by \cite{bo,gu}, we introduce the horizontal and vertical regularity and obtain the following theorem:
\begin{theorem}\label{t2.2}
Let $(\rho,\mathbf u)$ be a solution of CPE \eqref{1.1} in the sense of distribution.
Let us assume
$\rho,\mathbf v, \rho\mathbf v\in B^{\alpha,\infty}_3((0,T);B^{\alpha,\infty}_{3}(\mathbb{T};B^{\beta,\infty}_{3}(\mathbb{T}^2)))$,
where $\beta>\frac{2}{3}$, $\alpha>\frac{1}{3}$ and $2\alpha+\beta>2$.
Moreover, we assume further
\begin{align*}
0\leq\underline{\rho}\leq\rho\leq\overline{\rho},\hspace{3pt}{\text a.e. in}\hspace{2pt}(0,T)\times\mathbb{T}^3
\end{align*}
for some constants $\underline{\rho},\overline{\rho}$.
Moreover, we assume $w\in L^2(0,T;L^2(\mathbb{T}^3))$,
$\mathbf v\in L^{2}((0,T);H^1_{\bf x}(\mathbb{T}^2)L^2_z(\mathbb{T}^1))$ and $p\in C[\underline{\rho},\overline{\rho}]$.
Then the energy is locally conserved i.e.
\begin{align*}
\partial_t\left(\frac{1}{2}\rho|\mathbf v|^2+P(\rho)\right)+\rho|\nabla_{\bf x}\mathbf v|^2
+
\rho|\partial_z\mathbf v|^2
+
{\rm div}_{\bf x}\left[
\left(\frac{1}{2}\rho|\mathbf v|^2+p(\rho)+P(\rho)
-
\rho\nabla_{\bf x}\mathbf v\mathbb{I}
\right)\mathbf v\right]\\
+\partial_z\left[
\left(
\frac{1}{2}\rho|\mathbf v|^2+p(\rho)+P(\rho)
\right)
w
-
\rho\partial_z\mathbf v\cdot\mathbf v\right]=0,
\end{align*}
in the sense of distributions on $(0,T)\times\mathbb{T}^3$.
\end{theorem}

\begin{remark}

Notice that the assumptions in Theorem \ref{t2.2} about $\alpha$ and $\beta$ are the same as \cite{bo} due to the special structure of CPE. Moreover, from the assumptions $\beta>\alpha$ and $2\alpha+\beta>2$, we can deduce $\beta>\frac{2}{3}$.
It is important to point out that the procedure followed in the proof of Theorem \ref{t2.2} is different from \cite{bo} as our case is compressible and diffusive. 
\end{remark}

\begin{remark}
It should be noticed that vertical regularity of $\mathbf v\in L^{2}((0,T);H^1_{\bf x}(\mathbb{T}^2)L^2_z(\mathbb{T}^1))$
in the Theorem \ref{t2.2} is lower than  the regularity in the Theorem \ref{t2.1} and \cite{ak}, where $\mathbf v\in L^2((0,T);H^1(\mathbb{T}^3))$,
taking advantage of the anisotropic regularity and the special structure of CPE. 
\end{remark}

\begin{remark}
In both theorems, we prove the energy equality with the degenerate viscosity as $\mu(\rho)=\rho$.
It should be interesting to consider a general $\mu(\rho)$ to obtain the similar results. The difficulty is how to obtain the commutator estimates when mollifying the nonlinear viscosity $\mu(\rho)$.
\end{remark}


\section{Proof of Theorem \ref{t2.1} and Theorem \ref{t2.2}} \label{S3}

\subsection{Energy equality}
We follow the strategy from \cite{ak,e4}. Precisely speaking: first, we mollify the CPE in both space and time;
second, we derive the local energy equality for regularized quantities;
third, we estimate commutator errors generated by nonlinear terms;
and {\color{magenta} last}, passing $(\varepsilon,\delta) \rightarrow 0$ (see below), commutators tend to zero and the result is obtained in the original quantities.

Firstly, we mollify the momentum equation $\eqref{1.1}_2$ and obtain
\begin{align}\label{3.1}
\partial_t(\rho \mathbf{v})^\varepsilon+\textrm{div}_{\bf x}(\rho\mathbf{v}\otimes\mathbf{v})^\varepsilon
+\partial_z(\rho\mathbf vw)^\varepsilon+\nabla_{\bf x} p^\varepsilon(\rho)
=
{\rm div}_{\bf x}(\rho\nabla_{\bf x}\mathbf v)^\varepsilon+\partial_{z}(\rho\partial_z\mathbf v)^\varepsilon.
\end{align}
Notice that the pressure term in \eqref{3.1} does not belong to $C^2$, thus it does not satisfy the conditions in \cite{e4}.
Therefore, we need to mollify the pressure as follows (based on \cite{ak}):
Take a sequence $p^\delta\in C^2[\underline{\rho},\overline{\rho}]$ that converges uniformly to $p\in C[\underline{\rho},\overline{\rho}]$,
that is for each $\delta>0$
\begin{equation}\label{pres}
\|p^\delta-p\|_{L^\infty}\leq\delta.
\end{equation}
Replacing $p^\delta$ in \eqref{3.1}, we obtain
\begin{align}
\partial_t(\rho \mathbf{v})^\varepsilon&+\textrm{div}_{\bf x}(\rho\mathbf{v}\otimes\mathbf{v})^\varepsilon
+\partial_z(\rho\mathbf vw)^\varepsilon
+\nabla_{\bf x}
\left(
p^\delta(\rho)
\right)^\varepsilon\nonumber\\
&=
{\rm div}_{\bf x}(\rho\nabla_{\bf x}\mathbf v)^\varepsilon+\partial_{z}(\rho\partial_z\mathbf v)^\varepsilon
+\nabla_{\bf x}\Big{(}(p^\delta(\rho))^\varepsilon-p^\varepsilon(\rho)\Big{)}.
\nonumber
\end{align}
It is easy to obtain
\begin{align}\label{3.3}
\partial_t(\rho^\varepsilon \mathbf{v}^\varepsilon)
+
\textrm{div}_{\bf x}\big{(}(\rho \mathbf{v})^{\varepsilon} \otimes\mathbf{v}^\varepsilon\big{)}
&+
\partial_z\big{(}(\rho w)^\varepsilon\mathbf v^\varepsilon\big{)}
+
\nabla_{\bf x} p^\delta(\rho^\varepsilon)
-
{\rm div}_{\bf x}
\left(
\rho^\varepsilon\nabla_{\bf x}\mathbf v^\varepsilon
\right)
-
\partial_{z}
\left(
\rho^\varepsilon\partial_z\mathbf v^\varepsilon
\right)\nonumber\\
&
=
\partial_t
\Big(
\rho^\varepsilon \mathbf{v}^\varepsilon-(\rho \mathbf{v})^\varepsilon
\Big)
+
\textrm{div}_{\bf x}
\Big( (\rho \mathbf{v})^{\varepsilon} \otimes\mathbf{v}^\varepsilon
-
(\rho\mathbf{v}\otimes\mathbf{v})^\varepsilon
\Big)
\nonumber\\
&
\hspace{10pt}
+
\partial_z
\Big(
(\rho w)^\varepsilon\mathbf v^\varepsilon
-
(\rho\mathbf vw)^\varepsilon
\Big)
\nonumber\\
&
\hspace{10pt}
-{\rm div}_{\bf x}
\Big{(}
\rho^\varepsilon\nabla_{\bf x}\mathbf v^\varepsilon
-
(\rho\nabla_{\bf x}\mathbf v)^\varepsilon\big{)}
-
\partial_{z}\big{(}\rho^\varepsilon\partial_z\mathbf v^\varepsilon
-
(\rho\partial_z\mathbf v)^{\varepsilon}
\Big{)}
\nonumber\\
&
\hspace{10pt}
+
\nabla_{\bf x}
\Big{(} p^\delta(\rho^\varepsilon)-(p^\delta(\rho))^\varepsilon\Big{)}
+
\nabla_{\bf x}\Big{(}(p^\delta(\rho))^\varepsilon-p^\varepsilon(\rho)
\Big{)}
\nonumber\\
&=R^\varepsilon=\displaystyle\sum_{i=1}^7 R^\varepsilon_i.
\end{align}
Multiplying by $\mathbf v^\varepsilon$, we get
\begin{align}
\partial_t\rho^\epsilon|\mathbf v^\varepsilon|^2&
+
\rho^\varepsilon\partial_t\frac{|\mathbf v^\varepsilon|^2}{2}
+
{\rm div}_{\bf x}{ (\rho \mathbf{v})^{\varepsilon}} |\mathbf v^\varepsilon|^2
+
(\rho \mathbf{v})^{\varepsilon} \nabla_{\bf x}\frac{|\mathbf v^\varepsilon|^2}{2}
+
\partial_z(\rho w)^\varepsilon|\mathbf v^\varepsilon|^2
+
(\rho w)^\varepsilon\partial_z\frac{|\mathbf v^\varepsilon|^2}{2}
\nonumber\\
&+\rho^\varepsilon\mathbf v^\varepsilon\cdot\nabla_{\bf x}(P^\delta(\rho^\varepsilon))'
+\rho^\varepsilon|\nabla_{\bf x}\mathbf v^\varepsilon|^2+\rho^\varepsilon|\partial_z\mathbf v^\varepsilon|^2
-{\rm div}_{\bf x}(\rho^\varepsilon\nabla_{\bf x}\mathbf v^\varepsilon\cdot\mathbf v^\varepsilon)-\partial_z(\rho^\varepsilon\partial_z\mathbf v^\varepsilon\cdot\mathbf v^\varepsilon)
\nonumber\\
&=R^\varepsilon\mathbf v^\varepsilon \label{3.4}
\end{align}
where
\begin{align*}
\nabla_{\bf x}p^\delta(\rho^\varepsilon)=(p^\delta(\rho^\varepsilon))'\nabla_{\bf x}
\rho^\epsilon=\rho^\varepsilon(P^\delta(\rho^\varepsilon))''\nabla_{\bf x}\rho^\varepsilon,
=\rho^\varepsilon\nabla_{\bf x}(P^\delta(\rho^\varepsilon))'.
\end{align*}

Note that, considering vacuum, it is not possible to multiply (\ref{3.3}) by
 $\frac{(\rho \mathbf v)^\epsilon}{\rho^\epsilon}$ instead of $\mathbf v^\epsilon$ in order to obtain energy equality.

Mollifying the continuity equation $\eqref{1.1}_1$ in space and time, we have
\begin{align}\label{3.5}
\partial_t\rho^\varepsilon
+
{\rm div}_{\bf x}(\rho\mathbf v)^\varepsilon+\partial_z(\rho w)^\varepsilon=0.
\end{align}
It is easy to obtain the following equality
\begin{align}\label{3.6}
(\partial_t\rho^\varepsilon+{\rm div}_{\bf x}(\rho\mathbf v)^\varepsilon+\partial_z(\rho w)^\varepsilon)\frac{|\mathbf v^\varepsilon|^2}{2}=0.
\end{align}
We rewrite (\ref{3.5}) and get
\begin{align}
\partial_t\rho^\varepsilon+{\rm div}_{\bf x}(\rho^\varepsilon\mathbf v^\varepsilon)
+
\partial_z(\rho w)^\varepsilon
=
{\rm div}_{\bf x}(\rho^\varepsilon\mathbf v^\varepsilon-(\rho\mathbf v)^\varepsilon),
\nonumber\end{align}
then
\begin{align}\label{3.8}
\partial_tP^\delta(\rho^\varepsilon)
+
{\rm div}_{\bf x}(\rho^\varepsilon\mathbf v^\varepsilon)
(P^\delta(\rho^\varepsilon))'+\partial_z(\rho w)^\varepsilon(P^\delta(\rho^\varepsilon))'
={\rm div}_{\bf x}
\Big(
\rho^\varepsilon\mathbf v^\varepsilon-(\rho\mathbf v)^\varepsilon
\Big)(P^\delta(\rho^\varepsilon))'.
\end{align}

We define $s^\varepsilon={\rm div}_{\bf x}(\rho^\varepsilon\mathbf v^\varepsilon-(\rho\mathbf v)^\varepsilon)(P^\delta(\rho^\varepsilon))'$, and put \eqref{3.4}-\eqref{3.8} together to obtain
\begin{align}\label{3.9}
\partial_t
\left(
\frac{1}{2}\rho^\varepsilon|\mathbf v^\varepsilon|^2
+P^\delta(\rho^\varepsilon)
\right)
+&\rho^\varepsilon|\nabla_{\bf x}\mathbf v^\varepsilon|^2
+\rho^\varepsilon|\partial_z\mathbf v^\varepsilon|^2
\nonumber\\
+&
{\rm div}_{\bf x}
\left((
\rho\mathbf v)^{\varepsilon}
\frac{1}{2}|\mathbf v^\varepsilon|^2
+
(\rho^{\varepsilon} \mathbf v^\varepsilon) (P^\delta(\rho^\varepsilon))'\mathbf
-
\rho^\varepsilon\nabla_{\bf x}\mathbf v^\varepsilon\mathbf v^\varepsilon
\right)
\nonumber\\
+&\partial_z
\left(
(\rho w)^\varepsilon\frac{1}{2}|\mathbf v^\varepsilon|^2
+(\rho w)^\varepsilon(P^\delta(\rho^\varepsilon))'
-
\rho^\epsilon\partial_z\mathbf v^\epsilon\cdot\mathbf v^\varepsilon
\right)
 -
\left(R^\varepsilon\mathbf v^\varepsilon+s^\varepsilon\right)=0.
\end{align}
Therefore, we just need to prove the limit $\varepsilon\rightarrow0$ for each fixed $\delta>0$, then consider $\delta\rightarrow0$ in the following section.
\subsection{Commutator estimates}\label{ss32}
Section 3.1 can be seen as the prior estimates of energy equality. We will deduce the local energy equality for the regularized
quantities as follows:

Let $\varphi\in C^\infty_c(\Omega)$ be a test function and take $\varepsilon>0$ small enough so that $supp(\varphi)\subset(\varepsilon,T-\varepsilon)\times\mathbb{T}^3$.

Multiplying \eqref{3.3} with $\varphi\mathbf v^\varepsilon$ and integrating in time and space,
we have
\begin{align}
\int^\tau_0\int_{\mathbb{T}^3}\partial_t(\rho^\varepsilon\mathbf v^\varepsilon)\cdot\varphi\mathbf v^\varepsilon
\,
d{ \bf x}dzdt
&+\int^\tau_0\int_{\mathbb{T}^3}{\rm div}_{\bf x}\big{(}(\rho\mathbf{v})^\varepsilon\otimes\mathbf{v}^\varepsilon)\cdot\varphi\mathbf v^\varepsilon
+\partial_z\big{(}(\rho w)^\varepsilon\mathbf v^\varepsilon)\cdot\varphi\mathbf v^\epsilon
\,
d{ \bf x}dzdt
\nonumber
\\&
+\int^\tau_0\int_{\mathbb{T}^3}\nabla_{\bf x}p^\delta(\rho^\varepsilon)\cdot\varphi\mathbf v^\varepsilon
\,
d{ \bf x}dzdt \nonumber\\
&=\int^\tau_0\int_{\mathbb{T}^3}{\rm div}_{\bf x}(\rho^\varepsilon\nabla_{\bf x}\mathbf v^\varepsilon)\cdot\varphi\mathbf v^\varepsilon
+\partial_z(\rho^\varepsilon\partial_z\mathbf v^\varepsilon)\varphi\mathbf v^\varepsilon
\,
d{ \bf x}dzdt\nonumber\\
&+\displaystyle\sum_{i=1}^7\int^\tau_0\int_{\mathbb{T}^3}R^\varepsilon_i\cdot\varphi\mathbf v^\varepsilon
\,
d{ \bf x}dzdt.
\nonumber
\end{align}

Then by using the same argument as in \cite{ak,e4}, that is multiplying
(\ref{3.9}) by $\varphi$, and integrating in time and space, we obtain
\begin{align*}
\int^\tau_0\int_{\mathbb{T}^3}\partial_t
&\left(
\frac{1}{2}\rho^\varepsilon|\mathbf v^\varepsilon|^2
+P^\delta(\rho^\varepsilon)
\right)\varphi d{\bf x} dz dt
+\int^\tau_0\int_{\mathbb{T}^3}
\left(
\rho^\varepsilon|\nabla_{\bf x}\mathbf v^\varepsilon|^2
+\rho^\varepsilon|\partial_z\mathbf v^\varepsilon|^2
\right)\varphi d{\bf x} dz dt
\\
+&
\int^\tau_0\int_{\mathbb{T}^3}{\rm div}_{\bf x}
\left((
\rho\mathbf v)^{\varepsilon}
\frac{1}{2}|\mathbf v^\varepsilon|^2
+
(\rho^{\varepsilon} \mathbf v^\varepsilon) (P^\delta(\rho^\varepsilon))'\mathbf
-
\rho^\varepsilon (\nabla_{\bf x}\mathbf v^\varepsilon)\mathbf v^\varepsilon
\right)\varphi d{\bf x} dz dt
\\
+&\int^\tau_0\int_{\mathbb{T}^3}\partial_z
\left(
(\rho w)^\varepsilon\frac{1}{2}|\mathbf v^\varepsilon|^2
+(\rho w)^\varepsilon(P^\delta(\rho^\varepsilon))'
-
\int^\tau_0\int_{\mathbb{T}^3}\rho^\epsilon\partial_z\mathbf v^\epsilon\cdot\mathbf v^\varepsilon
\right)\varphi d{\bf x} dz dt\\
-&
\int^\tau_0\int_{\mathbb{T}^3}
\left(R^\varepsilon {\color{magenta} \cdot}\mathbf v^\varepsilon+s^\varepsilon\right)\varphi d{\bf x} dz dt=0.
\end{align*}

Therefore, we just need to show that each commutator error term converges to zero as
\begin{align*}
\int^\tau_0\int_{\mathbb{T}^3}R^\varepsilon_i {\color{magenta} \cdot}\mathbf v^\varepsilon\varphi \, d{ \bf x}dzdt
\rightarrow0 \quad
i=1,\dots,7
,\hspace{11pt}
\int^\tau_0\int_{\mathbb{T}^3}s^\varepsilon\varphi \, d{ \bf x}dzdt \rightarrow0.
\end{align*}
The following sections are devoted to prove Theorem \ref{t2.1} and Theorem \ref{t2.2}, respectively.}

\subsubsection{Proof of Theorem \ref{t2.1}}

The terms $R^\varepsilon_i$ $(i=1,2)$ and $s^\varepsilon$ are treated in the same way as in \cite{e4}.
For readers' convenience, we recall the main steps
for bounding $R^\varepsilon_2$ (the proof for $R^\varepsilon_1$ is similar, changing $\nabla_{\bf x}$ by $\partial_t$).

Using Remark \ref{r2.2}, we have
\begin{align}\label{3.11}
(\rho\mathbf v)^\varepsilon\otimes\mathbf v^\varepsilon
-(\rho\mathbf v\otimes\mathbf v)^\varepsilon
&=\big{(}(\rho\mathbf v)^\varepsilon-\rho\mathbf v\big{)}\otimes(\mathbf v^\varepsilon-\mathbf v)
\hfill
\nonumber\\
 &-\int^\varepsilon_{-\varepsilon}\int_{\mathbb{T}^3}
 \- \-
 \eta^\varepsilon(\tau,{\boldsymbol \xi_h},\xi_z)
\left(
\rho\mathbf v(t-\tau, {\bf x}- {\boldsymbol \xi}_h,z-\xi_z)
-\rho\mathbf v(t,{\bf x},z)
\right)
\nonumber
\\
&\hspace{2cm}
\otimes
\left(
\mathbf v(t-\tau, {\bf x}-{\boldsymbol \xi}_h, z-\xi_z)-\mathbf v(t, {\bf x},z)
\right) \,
d{\boldsymbol\xi}_h d\xi_z d\tau.
\end{align}
Recall that $R_2^{\varepsilon}=\textrm{div}_{\bf x}
\Big( (\rho \mathbf{v})^{\varepsilon} \otimes\mathbf{v}^\varepsilon
-
(\rho\mathbf{v}\otimes\mathbf{v})^\varepsilon
\Big)
$.
Therefore, for dealing with of $\int^\tau_0\int_{\mathbb{T}^3}R_2^{\varepsilon} {\color{magenta} \cdot}\mathbf v^{\varepsilon}\varphi \, d{ \bf x}dzdt
$,
we split it into two parts
($R_{2,1}^{\varepsilon}+R_{2,2}^{\varepsilon}$)
according to \eqref{3.11}.
\begin{align*}
R_{2,1}^{\varepsilon}&= \int^\tau_0\int_{\mathbb{T}^3}{\rm div}_{\bf x}\big{(}(\rho\mathbf v)^\varepsilon-\rho\mathbf v)\otimes(\mathbf v^\varepsilon-\mathbf v)\big{)}\varphi\mathbf v^\varepsilon d{ \bf x}dzdt\\
&\leq
\|\varphi\|_{C^1}\|(\rho\mathbf v)^\varepsilon-\rho\mathbf v\|_{L^3}\|\mathbf v^\varepsilon-\mathbf v\|_{L^3}\|\mathbf v^\varepsilon\|_{L^3}
+\|\varphi\|_{C^0}\|(\rho\mathbf v)^\varepsilon-\rho\mathbf v\|_{L^3}
\|\mathbf v^\varepsilon-\mathbf v\|_{L^3}\|\nabla_{\bf x}\mathbf v^\varepsilon\|_{L^3}\\
&\leq
C
\left(\varepsilon^{\alpha+\beta}
+
\varepsilon^{2\alpha+\beta-1}
\right)
\|\rho\mathbf v\|_{B_{3}^{\beta,\infty}}\|\mathbf v\|^2_{B_{3}^{\alpha,\infty}}.
\end{align*}
From the assumption in Theorem \ref{t2.1}, $\alpha>\frac{5}{6}$, it is easy to deduce $2\alpha+\beta>1$, so that we can obtain the term vanishes when $\varepsilon\rightarrow0$.
 For bounding $R_{2,2}^{\varepsilon}$, we use the Fubini's theorem and H\"{o}lder's inequality, and get
\begin{align}
R_{2,2}^{\varepsilon}&=
\int^\tau_0\int_{\mathbb{T}^3}
{\rm div}
\left(
\int^\varepsilon_{-\varepsilon}\int_{\mathbb{T}^3}\eta^\varepsilon(\tau,\xi)\right.
\Big(
\rho\mathbf v(t-\tau,{\bf x}-{\boldsymbol \xi}_h,z-\xi_z)-\rho\mathbf v(t,{\bf x},z)
\Big)
\nonumber
\\
&
\qquad \left.
\otimes\Big(
\mathbf v(t-\tau,{\bf x}-{\boldsymbol \xi}_h,z-\xi)-\mathbf v(t,{ \bf x},z)
\Big)
d{\boldsymbol \xi}_h \, d\xi_z \, d\tau
\right)  \cdot
\varphi\mathbf v^\varepsilon d{ \bf x}dzdt\nonumber\\
&
\leq
C
\left(
\varepsilon^{\alpha+\beta}
+
\varepsilon^{2\alpha+\beta-1}
\right)
\|\rho\mathbf v\|_{B_{3}^{\beta,\infty}}\|\mathbf v\|^2_{B_{3}^{\alpha,\infty}}.
\end{align}
Next, we turn to estimate the commutators $R^\varepsilon_6$ and $R^\varepsilon_7$.
The proof is inspired by \cite{ak} and we include the proof here for completeness.
\begin{align*}
R_{6}^{\varepsilon}&=\int^\tau_0\int_{\mathbb{T}^3}
\nabla_{\bf x}\big{(}p^\delta(\rho^\varepsilon)
-(p^\delta(\rho))^\varepsilon\big{)}\varphi\mathbf v^\varepsilon d{ \bf x}dzdt\\
&
=
-\int^\tau_0\int_{\mathbb{T}^3}\big{(}p^\delta(\rho^\varepsilon)
-
(p^\delta(\rho))^\varepsilon\big{)} (\varphi{\rm div}_{\bf x}\mathbf v^\varepsilon
+
\nabla_{\bf x}\varphi\mathbf v^\varepsilon)d{ \bf x}dzdt\\
&\leq
\|\varphi\|_{C^0}\|p^\delta(\rho^\varepsilon)-(p^\delta(\rho))^\varepsilon\|_{L^\frac{3}{2}}
\|{\rm div}_{\bf x}\mathbf v^\varepsilon\|_{L^{3}}
+
C\|\varphi\|_{C^1}\|p^\delta(\rho^\varepsilon)-(p^\delta(\rho))^\varepsilon\|_{L^\frac{3}{2}}
\|\mathbf v^\varepsilon\|_{L^3}\\
&\leq
C
\|p^\delta(\rho^\varepsilon)-(p^\delta(\rho))^\varepsilon\|_{L^\frac{3}{2}}
\varepsilon^{\alpha-1}
\|\mathbf v\|_{B^{\alpha,\infty}_3}
+
C
\|p^\delta(\rho^\varepsilon)-(p^\delta(\rho))^\varepsilon\|_{L^\frac{3}{2}}
\|\mathbf v^{\varepsilon}\|_{L^3}\\
&\leq
C
\big{(}\|\rho^\varepsilon-\rho\|^2_{L^{ 3}}
+
\sup_{{\bf y}\in{\rm supp}(\eta^\varepsilon)}\|\rho(\cdot)-\rho(\cdot-{\bf y})\|^2_{L^{ 3}}\big{)}
(\varepsilon^{\alpha-1} + 1)\|\mathbf v\|_{B^{\alpha,\infty}_3}
\\
& \le
C \, \varepsilon^{2\beta}
(\varepsilon^{\alpha-1} + 1)
  \|\rho\|^2_{B_3^{\beta,\infty}}
\|\mathbf v\|_{B^{\alpha,\infty}_3},
\end{align*}
where we have used Lemma \ref{l2.4a} to bound the terms on $p^\delta(\rho^\varepsilon)-(p^\delta(\rho))^\varepsilon$
and the fact that $p^{\delta}\in C^2[\underline{\rho},\overline{\rho}]$.
On the other hand, by using \eqref{pres}
\begin{align*}
R_{7}^{\varepsilon}&= \int^\tau_0\int_{\mathbb{T}^3}
\nabla_{\bf x}\big{(}(p^\delta(\rho))^\varepsilon-p^\varepsilon(\rho)\big{)}
\cdot \varphi \, \mathbf v^\varepsilon d{ \bf x} \, dz \, dt\\
&=-\int^\tau_0\int_{\mathbb{T}^3}
\big{(}p^\delta(\rho)-p(\rho)\big{)}^\varepsilon
\, \varphi
\, {\rm div}_{\bf x}\mathbf v^\varepsilon
d{ \bf x}dzdt
-\int^\tau_0\int_{\mathbb{T}^3}
\big{(}p^\delta(\rho)-p(\rho)
\big{)}^\varepsilon
\, \nabla_{\bf x}\varphi
\cdot
\mathbf v^\varepsilon
d{ \bf x} \, dz \, dt\\
&\leq
C \, \|\varphi\|_{C^0}\|\big{(}p^\delta(\rho)-p(\rho)\big{)}^\varepsilon\|_{L^\infty}
\|{\rm div}_{\bf x}\mathbf v^\varepsilon\|_{L^{2}}
+C\|\varphi\|_{C^1}\|\big{(}p^\delta(\rho)-p(\rho)\big{)}^\varepsilon\|_{L^\infty}\|\mathbf v^\varepsilon\|_{L^3}\\
&\leq C
\, \delta.
\end{align*}
The remaining is to estimate $R^\varepsilon_4$ and $R^\varepsilon_5$. We give a different proof compared with \cite{ak} as the following estimate for $R^\varepsilon_4$:
\begin{align}
R_{4}^{\varepsilon}&=\int^\tau_0
\int_{\mathbb{T}^3}{\rm div}_{\bf x}(\rho^\varepsilon\nabla_{\bf x}
\mathbf v^\varepsilon-(\rho\nabla_{\bf x}\mathbf v)^\varepsilon)
\cdot
\varphi\mathbf v^\varepsilon
d{ \bf x}dzdt\nonumber\\
&=- \int^\tau_0\int_{\mathbb{T}^3}(\rho^\varepsilon\nabla_{\bf x}
\mathbf v^\varepsilon-(\rho\nabla_{\bf x}\mathbf v)^\varepsilon)
\, : \, (\nabla_{\bf x}\varphi\mathbf v^\varepsilon )
d{ \bf x}dzdt
-\int^\tau_0\int_{\mathbb{T}^3}(\rho^\varepsilon\nabla_{\bf x}
\mathbf v^\varepsilon-(\rho\nabla_{\bf x}\mathbf v)^\varepsilon)
\, : \, \varphi\nabla_{\bf x}\mathbf v^\varepsilon
d{ \bf x}dzdt. \label{3.13b}
\end{align}
Using Remark \ref{r2.2}, we split the two terms on the right hand side of \eqref{3.13b} into four terms:
\begin{align*}
&\int^\tau_0\int_{\mathbb{T}^3}(\rho^\varepsilon\nabla_{\bf x}\mathbf v^\varepsilon-(\rho\nabla\mathbf v)^\varepsilon)
\, : \, (
\nabla_{\bf x}\varphi\mathbf v^\varepsilon) \,  d{ \bf x}dzdt\\
&=\int^\tau_0\int_{\mathbb{T}^3}(\rho^\varepsilon-\rho)
(\nabla_{\bf x}\mathbf v^\varepsilon-\nabla_{\bf x}\mathbf v)\mathbf v^\epsilon\cdot\nabla_{\bf x}\varphi
\, d{ \bf x}dzdt\\
&\hspace{3pt}
-\int^\varepsilon_{-\varepsilon}\int_{\mathbb{T}^3}
\- \- \eta^\varepsilon(\tau,\xi)
\Big{(}\rho(t-\tau, {\bf x}- {\boldsymbol \xi}_h)-\rho(t,{\bf x})
\Big{)}
\-
\Big{(}\nabla_{\bf x}\mathbf v(t-\tau,{\bf x}-{\boldsymbol \xi}_h,z-\xi_z)-\nabla_{\bf x}\mathbf v(t,x,z)\Big{)}\mathbf v^\epsilon\cdot\nabla_{\bf x}\varphi \,
d{\boldsymbol \xi}_h   d\xi_z  d\tau\\
&=J_1+J_2,
\end{align*}
and
\begin{align*}
&\int^\tau_0\int_{\mathbb{T}^3}(\rho^\varepsilon\nabla_{\bf x}\mathbf v^\varepsilon-(\rho\nabla_{\bf x}\mathbf v)^\varepsilon)\varphi\nabla_{\bf x}\mathbf v^\varepsilon d{ \bf x}dzdt\\
&=\int^\tau_0\int_{\mathbb{T}^3}(\rho^\varepsilon-\rho)(\nabla_{\bf x}\mathbf v^\varepsilon-\nabla_{\bf x}\mathbf v):\nabla_{\bf x}\mathbf v^\epsilon\varphi
\, d{ \bf x}dzdt\\
&\hspace{3pt}-\int^\varepsilon_{-\varepsilon}\int_{\mathbb{T}^3}\eta^\varepsilon(\tau,\xi)
\Big(
\rho(t-\tau,{\bf x}-{\boldsymbol \xi}_h)-\rho(t,{\bf x})
\Big)
\Big{(}\nabla_{\bf x}\mathbf v(t-\tau,{\bf x}-{\boldsymbol \xi}_h,z- \xi_z)-\nabla_{\bf x}\mathbf v(t,{\bf x},z)\Big{)}
\- : \- \nabla_{\bf x}\mathbf v^\epsilon\varphi \,
d{\boldsymbol \xi}_h  d{\xi_z}  d\tau\\
&=J_3+J_4.
\end{align*}
Observe that $J_2$ and $J_4$ has been written taking into account that $\rho$ is independent of $z$.

We just concern about the terms $J_1$ and $J_3$, the remaining ones can be estimated similarly by applying the Fubini's theorem (see in \cite{ak,e4}).
\begin{align*}
&|J_1|+|J_3|\\
&=
\left\vert
\int^\tau_0\int_{\mathbb{T}^3}(\rho^\varepsilon-\rho)(\nabla_{\bf x}\mathbf v^\varepsilon-\nabla_{\bf x}\mathbf v):\mathbf v^\epsilon\otimes\nabla_{\bf x}\varphi \, d{ \bf x}dzdt
\right\vert
+
\left\vert
\int^\tau_0\int_{\mathbb{T}^3}(\rho^\varepsilon-\rho)(\nabla_{\bf x}\mathbf v^\varepsilon-\nabla_{\bf x}\mathbf v):\nabla_{\bf x}\mathbf v^\epsilon\varphi \, d{ \bf x}dzdt
\right\vert\\
& \leq
\Vert \varphi \Vert_{C^1}
\Vert \rho^{\varepsilon}-\rho \Vert_{L^\infty}
\Vert \nabla_{\bf x}(\mathbf{v}^{\varepsilon}-\mathbf{v}) \Vert_{L^2}
\Vert \mathbf{v}^{\varepsilon} \Vert_{L^2}
+
\Vert \varphi \Vert_{C^0}
\Vert \rho^{\varepsilon}-\rho \Vert_{L^\infty}
\Vert \nabla_{\bf x}(\mathbf{v}^{\varepsilon}-\mathbf{v}) \Vert_{L^2}
\Vert \nabla_{\bf x} \mathbf{v}^{\varepsilon} \Vert_{L^2}.
\end{align*}

Recalling $\mathbf v\in L^2((0,T);H^1)$, we just use that $\Vert \nabla_{\bf x}(\mathbf{v}^{\varepsilon}-\mathbf{v}) \Vert_{L^2}\rightarrow0$ as $\varepsilon\rightarrow0$ and we deduce that $R^\varepsilon_4\rightarrow0$ as $\varepsilon\rightarrow0$. For the term $R^\varepsilon_5$, it is similar as $R^\varepsilon_4\rightarrow0$, hence we omit the detail.

Finally, we turn to analysis the difficult commutator error $R^\varepsilon_3$. As usual, the decomposition of $R^\varepsilon_3$ is similar to $R^\varepsilon_2$, so we divide it into two parts as the following:
\begin{align*}
&(\rho w)^\varepsilon\mathbf v^\varepsilon
 -(\rho w\mathbf v)^\varepsilon
=((\rho w)^\varepsilon-\rho w)(\mathbf v^\varepsilon-\mathbf v)\\
&-\int^\varepsilon_\varepsilon\int_{\mathbb{T}^3}\eta^\varepsilon(\tau,\xi)
\big(
\rho w(t-\tau, {\bf x}-{\boldsymbol \xi}_h,z-\xi_z)-\rho w(t,{ \bf x},z)
\big)
\big(\mathbf v(t-\tau,{\bf x}-{\boldsymbol \xi}_h, z-\xi_z)-\mathbf v(t,{\bf x},z)
\big)
\, d{\boldsymbol \xi}_h  d\xi_z   d\tau.
\end{align*}
We only focus on the first part, as the second part produces the same estimates after applying Fubini's theorem. It is difficult to estimate this term directly, because we do not have enough regularity of vertical velocity. This is the essential difficulty and distinction between Navier-Stokes system and CPE system.

Since we assume $0\leq\underline{\rho}\leq\rho\leq\overline{\rho}$ and $w\in L^{2}((0,T);L^{2}(\mathbb{T}^3))$, which means $\rho w\in L^{2}((0,T);L^{2}(\mathbb{T}^3))$, then ones get

\begin{align*}
\int^\tau_0\int_{\mathbb{T}^3}&((\rho w)^\varepsilon-\rho w)(\mathbf v^\varepsilon-\mathbf v)(\partial_z\varphi\mathbf v^\varepsilon+\varphi\partial_z\mathbf v^\varepsilon) \, d{ \bf x}dzdt\\
&\leq\|(\rho w)^\varepsilon-\rho w\|_{L^2}\|\mathbf v^\varepsilon-\mathbf v\|_{L^3}
\Big(
\|\varphi\|_{C^1}\|\mathbf v^\varepsilon\|_{L^6}+
\|\varphi\|_{C^0}\|\partial_z\mathbf v^\varepsilon\|_{L^6}
\Big)\\
&\leq C\|\rho w\|_{L^2}\varepsilon^{\alpha}\|\mathbf v\|_{B^{\alpha,\infty}_{3}}
\Big(
\|\mathbf v\|_{L^6}+\varepsilon^{\alpha-\frac{3}{2}}\|\mathbf v\|_{B^{\alpha-\frac{5}{3},\infty}_{6}}
\Big)\\
&\leq C
\left(
\varepsilon^{\alpha}
+
\varepsilon^{2\alpha-\frac{5}{3}}
\right)
\|\mathbf v\|^2_{B^{\alpha,\infty}_{3}},
\end{align*}
where we have used Lemma \ref{2.4} and Lemma \ref{l2.2} to get
\begin{align}\label{3.10}
\Vert \mathbf{v} \Vert_{B_6^{\alpha-2/3,\infty}} \le C \Vert \mathbf{v} \Vert_{B_3^{\alpha,\infty}},
\quad
\|\nabla_{\bf x} \mathbf v^\varepsilon\|_{L^6}\leq\varepsilon^{\alpha-\frac{5}{3}}\|\mathbf v\|_{B^{\alpha-\frac{1}{2},\infty}_{6}}
\leq \varepsilon^{\alpha-\frac{5}{3}}\|\mathbf v\|_{B^{\alpha,\infty}_{3}}.
\end{align}

\begin{remark}\label{r3.1}
In \eqref{3.10} we have used the third estimate in Lemma \ref{2.4} with $\omega $
$$\Vert\nabla_{\bf x} \mathbf v^\varepsilon\|_{L^6} \leq \varepsilon^{\omega - 1}\|\mathbf v\|_{B^{\omega,\infty}_{6}}.$$
 Setting $\omega =\alpha - 2/3$ we get \eqref{3.10}. Since $\alpha>\frac{5}{6}$, it converges to zero.
\end{remark}

\begin{remark}\label{t3.2}
Under the assumption $\mathbf v\in B^{\alpha,\infty}_4((0,T);B^{\alpha,\infty}_{4}(\mathbb{T}^3))$ (assumed in \cite{bo} for the incompressible and inviscid case), we can reprove the above term as
\begin{align*}
\int^\tau_0\int_{\mathbb{T}^3}&((\rho w)^\varepsilon-\rho w)(\mathbf v^\varepsilon-\mathbf v)(\partial_z\varphi\mathbf v^\varepsilon+\varphi\partial_z\mathbf v^\varepsilon) \, d{ \bf x}dzdt\\
&\leq\|(\rho w)^\varepsilon-\rho w\|_{L^2}\|\mathbf v^\varepsilon-\mathbf v\|_{L^4}
\Big(
\|\varphi\|_{C^1}\|\mathbf v^\varepsilon\|_{L^4}+
\|\varphi\|_{C^0}\|\partial_z\mathbf v^\varepsilon\|_{L^4}
\Big)\\
&\leq C\|\rho w\|_{L^2}\varepsilon^{\alpha}\|\mathbf v\|_{B^{\alpha,\infty}_{4}}
\Big(
\|\mathbf v\|_{L^4}+\varepsilon^{\alpha-1}\|\mathbf v\|_{B^{\alpha,\infty}_{4}}
\Big)\\
&\leq C
\left(
\varepsilon^{\alpha}
+
\varepsilon^{2\alpha-1}
\right)
\|\mathbf v\|^2_{B^{\alpha,\infty}_{4}},
\end{align*}
where we need $\alpha>\frac{1}{2}$.
\end{remark}


In order to complete the proof of Theorem \ref{t2.1}, we need to show the final level of convergence for $\delta$.
For each fixed $\delta>0$, we have the in the limit \eqref{3.9} as $\varepsilon\rightarrow0$
\begin{align}\label{3.13}
\partial_t
\left(
\frac{1}{2}\rho|\mathbf v|^2+P^\delta(\rho)
\right)
+\rho|\nabla_{\bf x}\mathbf v|^2+\rho|\partial_z\mathbf v|^2
+{\rm div}_{\bf x}
\left(
(\rho\mathbf v)\frac{1}{2}|\mathbf v|^2
+(p^\delta(\rho)
+P^\delta(\rho))\mathbf v
-\rho \, (\nabla_{\bf x}\mathbf v) \mathbf v
\right)
\nonumber\\
+\partial_z
\left(
(\rho w)\frac{1}{2}|\mathbf v|^2
+(p^\delta(\rho)+P^\delta(\rho))w
-\rho\partial_z\mathbf v\cdot\mathbf v
\right)=0,
\end{align}
where $P^\delta(\rho)=\rho\displaystyle\int^\rho_1\frac{p^\delta(r)}{r^2}dr$ is the pressure potential.

\bigskip

We will prove that (\ref{3.13}) converges in the sense of distributions on $\Omega$ as $\delta\rightarrow0$ to \eqref{energy1}. This limit process of $\delta$ is similar to  Akramov et al. \cite{ak}. Taking $\varphi\in C^\infty_c(\Omega)$, we have
\begin{align*}
\Big\vert
\int^\tau_0\int_{\mathbb{T}^3}
(p^\delta(\rho)-p(\rho))
({\bf v} \cdot
\nabla_{\bf x}\varphi )
\,d{ \bf x}dzdt
\Big\vert
\leq
C \,
\|\varphi\|_{C^1}\|p^\delta(\rho)-p(\rho)\|_{L^\infty}\| {\bf v}\|_{L^3}\\
\leq
C \,
\|\varphi\|_{C^1}\|p^\delta(\rho)-p(\rho)\|_{L^\infty}\| {\bf v} \|_{B_3^{\alpha,\infty}}\leq C\delta.
\end{align*}

Recalling the estimates from \cite{ak}
\begin{align*}
|P^\delta(\rho)-P(\rho)|\leq\rho\int^\rho_1\frac{|p^\delta(r)-p(r)|}{r^2}dr\leq
\|p^\delta-p\|_{L^\infty}\rho
\left\vert
\int^\rho_1\frac{1}{r^2}dr
\right\vert
\leq(1+\rho)\|p^\delta-p\|_{L^\infty},
\end{align*}
we can obtain that

\begin{align*}
\Big\vert
\int^\tau_0\int_{\mathbb{T}^3}
(P^\delta(\rho)-P(\rho))
(\partial_{\bf x} \varphi \cdot {\bf v})
\, d{ \bf x}dzdt
\Big\vert
\leq\|\varphi\|_{C^1}
\Big(
1+\|\rho\|_{L^2}
\Big)
\|p^\delta-p\|_{L^\infty}\| {\bf v}\|_{L^2}\leq C\delta.
\end{align*}
The term in $\partial_z$ in (\ref{3.13}) can be estimated similarly. Then taking $\epsilon\rightarrow0$, we find that all commutator error terms converge to zero, which complete the proof of Theorem \ref{t2.1}.


\subsubsection{Proof for Theorem \ref{t2.2}}

Comparing with Theorem \ref{t2.1}, we delicately combine the anisotropic regularity with the special structure of CPE system, and observe some new results.
The process of proving Theorem \ref{t2.2} is similar to Theorem \ref{t2.1}.
We just focus on some different
estimates for commutators  $R^\varepsilon_2$ and $R^\varepsilon_6$.
First, we observe for $R^\varepsilon_2$
\begin{align*}
\int^\tau_0\int_{\mathbb{T}^3}&{\rm div}_{\bf x}\big{(}(\rho\mathbf v)^\varepsilon-\rho\mathbf v)\otimes(\mathbf v^\varepsilon-\mathbf v)\big{)}\varphi\mathbf v^\varepsilon \,
d{ \bf x}dzdt\\
&\leq
\|\varphi\|_{C^1}\|(\rho\mathbf v)^\varepsilon-\rho\mathbf v\|_{L^3}\|\mathbf v^\varepsilon-\mathbf v\|_{L^3}\|\mathbf v^\varepsilon\|_{L^3}
+\|\varphi\|_{C^0}\|(\rho\mathbf v)^\varepsilon-\rho\mathbf v\|_{L^3}
\|\mathbf v^\varepsilon-\mathbf v\|_{L^3}\|\nabla_{\bf x}\mathbf v^\varepsilon\|_{L^3}\\
&\leq
C\varepsilon^{2\alpha}
\|\rho\mathbf v\|_{B_{3,t}^{\alpha}(B_{3,z}^{\alpha}(B_{3,h}^{\beta}))}\|\mathbf v\|^2_{B_{3,t}^{\alpha}(B_{3,z}^{\alpha}(B_{3,h}^{\beta}))}\\
&\hspace{10pt}
+C\varepsilon^{2\alpha}
\|\rho\mathbf v\|_{B_{3,t}^{\alpha}(B_{3,z}^{\alpha}(B_{3,h}^{\beta}))}\|\mathbf v\|_{B_{3,t}^{\alpha},B_{3,z}^{\alpha}(B_{3,h}^{\beta})}
\varepsilon^{\alpha-1}
\|\mathbf v\|_{B_{3,t}^{\alpha}(B_{3,z}^{\alpha}(B_{3,h}^{\beta}))}
\\
&\leq
C
\left(
\varepsilon^{2\alpha}+\varepsilon^{3\alpha-1}
\right)
\|\rho\mathbf v\|_{B_{3,t}^{\alpha}(B_{3,z}^{\alpha}(B_{3,h}^{\beta}))}\|\mathbf v\|^2_{B_{3,t}^{\alpha}(B_{3,z}^{\alpha}(B_{3,h}^{\beta}))}.
\end{align*}
Here and in the following, we use the notation $B_{3,t}^{\alpha}(B_{3,z}^{\alpha}(B_{3,h}^{\beta}))$ instead of $ B^{\alpha,\infty}_3((0,T);B^{\alpha,\infty}_{3}(\mathbb{T};B^{\beta,\infty}_{3}(\mathbb{T}^2)))$.

Regarding the commutator term $R^\varepsilon_4$, we divide it into two parts as
\begin{align*}
\int^\tau_0
&\int_{\mathbb{T}^3}{\rm div}_{\bf x}(\rho^\varepsilon\nabla_{\bf x}
\mathbf v^\varepsilon-(\rho\nabla_{\bf x}\mathbf v)^\varepsilon)\varphi\mathbf v^\varepsilon
d{ \bf x}dzdt\nonumber\\
&=- \int^\tau_0\int_{\mathbb{T}^3}(\rho^\varepsilon\nabla_{\bf x}
\mathbf v^\varepsilon-(\rho\nabla_{\bf x}\mathbf v)^\varepsilon)\nabla_{\bf x}\varphi\mathbf v^\varepsilon d{ \bf x}dzdt
-\int^\tau_0\int_{\mathbb{T}^3}(\rho^\varepsilon\nabla_{\bf x}
\mathbf v^\varepsilon-(\rho\nabla_{\bf x}\mathbf v)^\varepsilon)\varphi\nabla_{\bf x}\mathbf v^\varepsilon
d{ \bf x}dzdt.
\end{align*}
Then the following process is the same as in the Theorem 2.1, so we just concern on the terms $J_1$ and $J_3$, that is
\begin{align*}
&|J_1|+|J_3|\\
&=
\left\vert
\int^\tau_0\int_{\mathbb{T}^3}(\rho^\varepsilon-\rho)(\nabla_{\bf x}\mathbf v^\varepsilon-\nabla_{\bf x}\mathbf v):\mathbf v^\epsilon\otimes\nabla_{\bf x}\varphi \, d{ \bf x}dzdt
\right\vert
+
\left\vert
\int^\tau_0\int_{\mathbb{T}^3}(\rho^\varepsilon-\rho)(\nabla_{\bf x}\mathbf v^\varepsilon-\nabla_{\bf x}\mathbf v):\nabla_{\bf x}\mathbf v^\epsilon\varphi \, d{ \bf x}dzdt
\right\vert\\
& \leq
\Vert \varphi \Vert_{C^1}
\Vert \rho^{\varepsilon}-\rho \Vert_{L^\infty}
\Vert \nabla_{\bf x}(\mathbf{v}^{\varepsilon}-\mathbf{v}) \Vert_{L^2}
\Vert \mathbf{v}^{\varepsilon} \Vert_{L^2}
+
\Vert \varphi \Vert_{C^0}
\Vert \rho^{\varepsilon}-\rho \Vert_{L^\infty}
\Vert \nabla_{\bf x}(\mathbf{v}^{\varepsilon}-\mathbf{v}) \Vert_{L^2}
\Vert \nabla_{\bf x} \mathbf{v}^{\varepsilon} \Vert_{L^2}.
\end{align*}

Recalling $\mathbf v\in L^2((0,T);H^1_{\mathbf x}(\mathbb{T}^2)L^2_z(\mathbb{T}))$ and the definition of  \eqref{2.1}, we just use that $\Vert \nabla_{\bf x}(\mathbf{v}^{\varepsilon}-\mathbf{v}) \Vert_{L^2}\rightarrow0$ as $\varepsilon\rightarrow0$,
then deduce that $R^\varepsilon_4\rightarrow0$ as $\varepsilon\rightarrow0$.

For the commutator term $R^\varepsilon_5$, that contains vertical derivatives, we use the anisotropic regularity to give a different and interesting proof. Recalling the density is independent of $z$, we integrate $R^\varepsilon_5$ by parts and rewrite it as
\begin{align*}
\int^\tau_0&
\int_{\mathbb{T}^3}\partial_z
\Big(
\rho^\varepsilon\partial_z\mathbf v^\varepsilon-(\rho\partial_z\mathbf v)^\varepsilon
\Big)
\cdot
\varphi\mathbf v^\varepsilon \-
d{ \bf x}dzdt\\
&=
\int^\tau_0\int_{\mathbb{T}^3}\partial_z
\Big(
\partial_z(\rho^\varepsilon\mathbf v^\varepsilon)-\partial_z(\rho\mathbf v)^\varepsilon
\Big)
\cdot
\varphi\mathbf v^\varepsilon \-
d{ \bf x}dzdt\\
&=-\int^\tau_0\int_{\mathbb{T}^3}
\Big(
\partial_z(\rho^\varepsilon\mathbf v^\varepsilon)-\partial_z(\rho\mathbf v)^\varepsilon
\Big)
\cdot
\Big(
\mathbf v^\epsilon\partial_z\varphi +\partial_z\mathbf v^\epsilon\varphi
\Big)
\, d{ \bf x}dzdt\\
&=-\int^\tau_0\int_{\mathbb{T}^3}
\Big(
\partial_z(\rho\mathbf v^\varepsilon)-\partial_z(\rho\mathbf v)^\varepsilon
\Big)
\cdot
\Big(
\mathbf v^\epsilon\partial_z\varphi +\partial_z\mathbf v^\epsilon\varphi
\Big)
\, d{ \bf x}dzdt\\
&\hspace{5pt}-\int^\tau_0\int_{\mathbb{T}^3}
\partial_z\Big(
(\rho^\varepsilon-\rho)\mathbf v^\varepsilon
\Big)
\cdot
\Big(
\mathbf v^\epsilon\partial_z\varphi +\partial_z\mathbf v^\epsilon\varphi
\Big)
 \, d{ \bf x}dzdt.
\end{align*}

Due to the independence of $\rho$ with respect to $z$, Besov estimates for $\rho$
only depends on variations in the $\beta$-index, which is the index related to the horizontal components. Based on this fact, we bound the second term as the following
\begin{align*}
\Big\vert
\int^\tau_0&\int_{\mathbb{T}^3}
\partial_z\Big(
(\rho^\varepsilon-\rho)\mathbf v^\varepsilon
\Big)
\cdot
\Big(
\mathbf v^\epsilon\partial_z\varphi +\partial_z\mathbf v^\epsilon\varphi
\Big)
\, d{ \bf x}dzdt
\Big\vert\\
&=
\Big\vert\int^\tau_0\int_{\mathbb{T}^3}
\Big(
(\rho^\varepsilon-\rho)(\partial_z\mathbf v^\varepsilon)
\Big)
\cdot
\Big(
\mathbf v^\epsilon\partial_z\varphi +\partial_z\mathbf v^\epsilon\varphi
\Big)
\, d{ \bf x}dzdt
\Big\vert\\
&\leq C\|\rho^\epsilon-\rho\|_{L^3}\|\partial_z\mathbf v^\epsilon\|_{L^3}
\Big(
\|\mathbf v^\epsilon\|_{L^3}+\|\partial_z\mathbf v^\epsilon\|_{L^3}
\Big)
\\
&\leq C
\varepsilon^{\beta}
\|\rho\|_{B^{\alpha}_{3,t}(B^\alpha_{3,z}(B^\beta_{3,h}))}
\varepsilon^{\alpha-1}\|\mathbf v\|_{B^{\alpha}_{3,t}(B^\alpha_{3,z}(B^\beta_{3,h}))}
\left(
\|\mathbf v\|_{L^3}+
\varepsilon^{\alpha-1}\|\mathbf v\|_{B^{\alpha}_{3,t}(B^{\alpha}_{3,z}(B^{\beta}_{3,h}))}
\right)\\
&\leq C
\Big(
\varepsilon^{\alpha+\beta-1} + \varepsilon^{2\alpha+\beta-2}
\Big)
\|\rho\|_{B^{\alpha}_{3,t}(B^\alpha_{3,z}(B^\beta_{3,h}))}
\|\mathbf v\|^2_{B^{\alpha}_{3,t}(B^\alpha_{3,z}(B^\beta_{3,h}))}.
\end{align*}

For the first term, we use H\"{o}lder's inequality to obtain
\begin{align*}
\Big\vert
\int^\tau_0&\int_{\mathbb{T}^3}
\Big(
\partial_z(\rho\mathbf v^\varepsilon)-\partial_z(\rho\mathbf v)^\varepsilon
\Big)
\cdot
\Big(
\mathbf v^\epsilon\partial_z\varphi +\partial_z\mathbf v^\epsilon\varphi
\Big)
\,
d{ \bf x}dzdt
\Big\vert\\
&\leq \|\partial_z(\rho\mathbf v^\varepsilon)-\partial_z(\rho\mathbf v)^\varepsilon\|_{L^3}
\Big(
\|\mathbf v^\epsilon\|_{L^3}\|\partial_z\varphi\|_{L^3}+\|\partial_z\mathbf v^\epsilon\|_{L^3}\|\varphi\|_{L^3}
\Big).
\end{align*}
Recalling Lemma \ref{l2.5}, we have $\|\partial_z(\rho\mathbf v^\varepsilon)-\partial_z(\rho\mathbf v)^\varepsilon\|_{L^3}\rightarrow0$ as $\varepsilon\rightarrow0$, where we take $\rho=f$ and $\mathbf v=\mathbf g$.

Then for $R^\varepsilon_6$, it is similar to obtain
\begin{align*}
\Big\vert
\int^\tau_0\int_{\mathbb{T}^3}&\nabla_{\bf x}\big{(}p^\delta(\rho^\varepsilon)
-(p^\delta(\rho))^\varepsilon\big{)}\varphi\mathbf v^\varepsilon
\, d{ \bf x}dzdt
\Big\vert
\\
&
=
\Big\vert
-\int^\tau_0\int_{\mathbb{T}^3}\big{(}p^\delta(\rho^\varepsilon)
-
(p^\delta(\rho))^\varepsilon\big{)} (\varphi{\rm div}_{\bf x}\mathbf v^\varepsilon
+
\nabla_{\bf x}\varphi\mathbf v^\varepsilon)
\, d{ \bf x}dzdt
\Big\vert
\\
&\leq
\|\varphi\|_{C^0}\|p^\delta(\rho^\varepsilon)-(p^\delta(\rho))^\varepsilon\|_{L^\frac{3}{2}}
\|{\rm div}_{\bf x}\mathbf v^\varepsilon\|_{L^3}
+
C\|\varphi\|_{C^1}\|p^\delta(\rho^\varepsilon)-(p^\delta(\rho))^\varepsilon\|_{L^\frac{3}{2}}
\|\mathbf v^\varepsilon\|_{L^3}\\
&\leq
C
\|p^\delta(\rho^\varepsilon)-(p^\delta(\rho))^\varepsilon\|_{L^{3/2}}
\varepsilon^{\alpha-1}
\|\mathbf v^{\varepsilon}\|_{B^{\alpha}_{3,t}(B^{\alpha}_{3,z}(B^{\beta}_{3,h}))}
+
C
\|p^\delta(\rho^\varepsilon)-(p^\delta(\rho))^\varepsilon\|_{L^{3/2}}
\|\mathbf v\|_{L^3}\\
&\leq
C
\big{(}\|\rho^\varepsilon-\rho\|^2_{L^3}
+
\sup_{y\in{\rm supp}\eta^\varepsilon}\|\rho(\cdot)-\rho(\cdot-y)\|^2_{L^3}\big{)}
\Big(
\varepsilon^{\alpha-1} + 1
\Big)
\|\mathbf v\|_{B^{\alpha}_{3,t}(B^\alpha_{3,z}(B^\beta_{3,h}))}
\\
&\leq
C  \varepsilon^{2\alpha}
\Big(
\varepsilon^{\alpha-1} + 1
\Big)
  \|\rho\|^2_{B^{\alpha}_{3,t}(B^\alpha_{3,z}(B^\beta_{3,h}))}
\|\mathbf v\|_{B^{\alpha}_{3,t}(B^\alpha_{3,z}(B^\beta_{3,h}))}      \\
&\leq
C \Big(
\varepsilon^{3\alpha-1} + \varepsilon^{2\alpha}
\Big)
  \|\rho\|^2_{B^{\alpha}_{3,t}(B^{\alpha}_{3,z}(B^{\beta}_{3,h}))}
\|\mathbf v\|_{B^{\alpha}_{3,t}(B^{\alpha}_{3,z}(B^{\beta}_{3,h}))}.
\end{align*}
Recalling the assumption $\alpha>\frac{1}{3}$, we could deduce this term converges to zero.

Let us stress that for the term $R^\varepsilon_3$, we need to apply the important property of (CPE). The following equality plays a key role in the existence of CPE \cite{liu2}, which is also helpful and important in our analysis,
\begin{align}
(\rho w)
{\color{blue}
(\x,z,t)
}
&=-\int^z_0{\rm div}_{\bf x}[\rho(\mathbf v(\mathbf x,z',t)-\int^1_0\mathbf v(\mathbf x,s,t)ds)]dz'. \nonumber\\
&=-\int^z_0{\rm div}_{\bf x}(\rho\widetilde{\mathbf v}) {\color{blue}
(\x,z',t)
} \, dz', \label{divPE}
\end{align}
where $\widetilde{\mathbf v}=\mathbf v-\overline{\mathbf v}=v-\int^1_0\mathbf v(x,s,t)ds$. We give the detail of the derivation of the above equation in the Appendix.

We will divide the proof into two cases according to the range of $\beta$.
When $\frac{2}{3}<\beta<1$, followed the same process of Theorem \ref{t2.1} to control the convergence of $R^\varepsilon_3$, we could obtain the corresponding result with $\alpha>\frac{5}{6}$. If $\beta>1$, due to the assumption $\rho\mathbf v\in B^{\alpha,\infty}_3((0,T);B^{\alpha,\infty}_{3}(\mathbb{T};B^{\beta,\infty}_{3}(\mathbb{T}^2)))$, then we can deduce $\rho w\in B^{\alpha,\infty}_{3}$$((0,T);B^{\alpha+1,\infty}_{3}(\mathbb{T};B^{\beta-1,\infty}_{3}(\mathbb{T}^2))$. It means $\rho w$ has $B^{\alpha+1,\infty}_{3}$ regularity in the $z-$ direction, $B^{\beta-1,\infty}_{3}$ regularity in the horizontal directions, and $B^{\alpha,\infty}_{3}$ regularity in time. Therefore, we can estimate $R^\varepsilon_3$ as follows
\begin{align*}
\int^\tau_0\int_{\mathbb{T}^3}&
\left(
(\rho w)^\varepsilon-\rho w)
(\mathbf v^\varepsilon-\mathbf v)
(\partial_z\varphi\mathbf v^\varepsilon+\varphi\partial_z\mathbf v^\varepsilon
\right)
\, d{ \bf x}dzdt\\
&\leq
\left(
\int_0^{\tau}\|(\rho w)^\varepsilon-\rho w\|^3_{L^3(\mathbb{T}^3)}dt
\right)^{\frac{1}{3}}
\cdot\|\mathbf v^\varepsilon-\mathbf v\|_{L^3}
\Big(
\|\varphi\|_{C^1}\|\|\mathbf v^\varepsilon\|_{L^3}+
\|\varphi\|_{C^0}\|\|\partial_z\mathbf v^\varepsilon\|_{L^3}
\Big)
\\
&\leq C
\varepsilon^{\beta-1}
\left(
\int_0^{\tau}\|\rho w\|^3_{B^{\alpha+1}_{3,z}(B^{\beta-1}_{3,h})}dt
\right)^{\frac{1}{3}}
\varepsilon^{\alpha}\|\mathbf v\|^2_{B^{\alpha}_{3,t}(B^{\alpha}_{3,z}(B^{\beta}_{3,h}))}\\
&\hspace{10pt}+C\varepsilon^{\beta-1}
\left(
\int_0^{\tau}\|\rho w\|^3_{B^{\alpha+1}_{3,z}(B^{\beta-1}_{3,h})}dt
\right)^{\frac{1}{3}}
\varepsilon^{\alpha}\|\mathbf v\|_{B^{\alpha}_{3,t}(B^{\alpha}_{3,z}(B^{\beta}_{3,h}))}
\varepsilon^{\alpha-1}\|\mathbf v\|_{ B^{\alpha}_{3,t}(B^{\alpha}_{3,z}(B^{\beta}_{3,h}))}\\
&=C
\left(
\varepsilon^{\alpha+\beta-1}
+
\varepsilon^{2\alpha+\beta-2}
\right)
\Big(
\int_0^{\tau}\|\rho w\|^3_{B^{\alpha+1}_{3,z}(B^{\beta-1}_{3,h})}dt
\Big)^{\frac{1}{3}}
\|\mathbf v\|^2_{B^{\alpha}_{3,t}(B^{\alpha}_{3,z}(B^{\beta}_{3,h}))}\\
&=C
\Big(
\varepsilon^{\alpha+\beta-1}
+
\varepsilon^{2\alpha+\beta-2}
\Big)
\|\rho w\|_{B^{\alpha}_{3,t}(B^{\alpha+1}_{3,z}(B^{\beta-1}_{3,h}))}
\|\mathbf v\|^2_{B^{\alpha}_{3,t}(B^{\alpha}_{3,z}(B^{\beta}_{3,h}))}.
\end{align*}
Then taking $\epsilon\rightarrow0$, we finish the proof of Theorem \ref{t2.2}.


\section{Conclusion}

Theorem \ref{t2.1} can be seen as a kind of generalization of Feireisl et al.~\cite{e4} and  Wiedemann et al.~\cite{ak} results for compressible Euler and Navier-Stokes equations, respectively, to the CPE case by using different proofs of some part of estimates (like $R^\varepsilon_4$,
$R^\varepsilon_5$) together with the observation of the special structure of CPE.
The classical isotropic case is analyzed in \cite{ak,e4}. But the conventional isotropic spaces do not reflect the diverse regularity in the horizontal and the vertical directions of the velocity field for CPE model.
And we should stress that since the regularity of $w$ is missing, which usually follows from the momentum equation of Navier-Stokes or Euler system, the only way to obtain the information of $w$ is through the continuity equation (see (\ref{divPE})). Comparing with Theorem \ref{t2.1}, Theorem \ref{t2.2} also considers the continuity of the pressure and anisotropic spaces in the horizontal and vertical directions for the density, horizontal velocity and horizontal momentum. Such anisotropic spaces already appear in Titi et al.~\cite{bo}, but, the CPE model we consider in this manuscript contains dissipation terms, vacuum and it is compressible, in contrast with \cite{bo} that treats the incompressible, inviscid and no vacuum case.


\section{Appendix}

This appendix is devoted to explore the relation between $\rho w$ and $\rho {\bf v}$, thanks to the continuity equation (\ref{1.1})$_1$.

\subsection{Expression of  $\rho w$ in terms of $\rho {\bf v}$}\label{A1}
We start from the continuity equation
\begin{equation}\label{f1}
\partial_t \rho ({\bf x},t)
+
{\rm div}_{\x} \left[
\rho(\x,t) \, \vi(\x,z,t)
\right]
+
\partial_z
\left[
\rho(\x,t) \, w(\x, z,t)
\right] =0,
\end{equation}
where we have taken into account that $\rho$ is independent of $z$.
Observe that such fact can be deduced from
\begin{equation}\label{f2}
\partial_z p(\rho)=\partial_z (\rho^{\gamma})=\gamma \, \rho^{\gamma-1} \partial_z \rho =0.
\end{equation}

\bigskip

Integrating (\ref{f1}) in $z \in (0,h)$, being $(0,h)$ the vertical ``reference" length of $\mathbb{T}^3$, we get:
\begin{equation}\label{f3}
h \, \partial_t \rho (\x,t)
+
\displaystyle\int_0^h {\rm div}_{\x} \left[
\rho(\x,t) \, \vi(\x,z,t)
\right] \, dz
+
\left[
\rho(\x,t) \, w(\x, h,t)
-
\rho(\x,t) \, w(\x, 0,t)
\right] =0.
\end{equation}
Taking into account the boundary conditions considered (remember that $\Omega=\mathbb{T}^3$, and therefore $w(\x, h,t)
=w(\x, 0,t)
$), then:
\begin{equation}\label{f4}
h \, \partial_t \rho (\x,t)
+
\displaystyle\int_0^h {\rm div}_{\x} \left[
\rho(\x,t) \, \vi(\x,z,t)
\right] \, dz
=0,
\end{equation}
which can be rewritten as:
\begin{equation}\label{f4-bis}
h \, \partial_t \rho (\x,t)
+
 {\rm div}_{\x} \left[
\rho(\x,t) \, \displaystyle\int_0^h \vi(\x,z,t) \, dz
\right] =0,
\end{equation}
and therefore
\begin{equation}\label{f4-tris}
\partial_t \rho (\x,t)
=-\displaystyle\frac{1}{h} \,
 {\rm div}_{\x} \left[
\rho(\x,t) \, \displaystyle\int_0^h \vi(\x,z,t) \, dz
\right]
=
-{\rm div}_{\x} \left[
\rho(\x,t) \, \bar{\vi}(\x,t)
\right]
\end{equation}
for
$$
\bar{\vi}(\x,t) = \displaystyle\frac{1}{h}
\,  \displaystyle\int_0^h \vi(\x,z,t) \, dz
$$

\bigskip

Now, we integrate (\ref{f1}) with respect to the vertical coordinate in $(0,z)$:
\begin{equation}\label{f5}
\displaystyle\int_0^z \partial_t \rho (\x,t)  \, dz
+
\displaystyle\int_0^z {\rm div}_{\x} \left[
\rho(\x,t) \, \vi(\x,s,t)
\right] \, ds
+
\left[
\rho(\x,t) \, w(\x, z,t)-
\rho(\x,t) \, w(\x, 0,t)
\right] =0.
\end{equation}
By using (\ref{f4-tris}), we rewrite (\ref{f5}) as:
\begin{equation}\label{f6}
\displaystyle\int_0^z {\rm div}_{\x} \left[
\rho(\x,t) \,
\left(
\vi(\x,s,t)
-
\bar{\vi}(\x,t)
\right)
\right] \, ds
+
\left[
\rho(\x,t) \, w(\x, z,t)-
\rho(\x,t) \, w(\x, 0,t)
\right] =0.
\end{equation}

\bigskip

Therefore,
\begin{equation}\label{f7}
\rho(\x,t) \,  w(\x, z,t) =\rho(\x,t)w(\x, 0,t)
-
\displaystyle\int_0^z {\rm div}_{\x} \left[
\rho(\x,t) \,
\left(
\vi(\x,s,t)
-
\bar{\vi}(\x,t)
\right)
\right] \, ds.
\end{equation}

Here we should emphasize that the physical boundary condition $w|_{z=0}=0$ is important to derive the formulation \eqref{divPE}. However, Titi et al. \cite{bo,liu3} pointed out that if the weak solution for PE and CPE satisfy suitable regularity assumptions and subject to some symmetry properties in the periodic boundary conditions, then the solutions automatically satisfy $w|_{z=0}=0$ when we restrict the solutions to the domain $\mathbb{T}^2\times(0,1)$. Therefore, we can deduce that
\begin{equation}\label{f8}
\rho(\x,t)\, w(\x, z,t) =-
  \,
\displaystyle\int_0^z {\rm div}_{\x} \left[
\rho(\x,t) \,
\left(
\vi(\x,s,t)
-
\bar{\vi}(\x,t)
\right)
\right] \, ds.
\end{equation}


\vskip 0.5cm


\vskip 0.5cm
\noindent {\bf Acknowledgements}

\vskip 0.1cm
The research of  \v{S}. Ne\v{c}asov\'{a} has been supported by the Czech Science Foundation (GA\v CR) project 22-01591S.  Moreover,  \v S. N.  has been supported by  Praemium Academiae of \v S. Ne\v casov\' a. The Institute of Mathematics, CAS is supported by RVO:67985840.
M.A. Rodr\'{\i}guez-Bellido has been supported by Grants PID2023-149182NB-I00 funded by MICIU/AEI/10.13039/501100011033 and by FEDER, UE, and
Grant US-1381261 (US/JUNTA/FEDER, UE).
 T. Tang is partially supported by NSFC No. 12371246 and NSF of Jiangsu Province Grant No. BK20221369.


\end{document}